# TESTING STATISTICAL HYPOTHESIS ON RANDOM TREES AND APPLICATIONS TO THE PROTEIN CLASSIFICATION PROBLEM[1]

By Jorge R. Busch, Pablo A. Ferrari, Ana Georgina Flesia,[2]
Ricardo Fraiman, Sebastian P. Grynberg and
Florencia Leonardi[3]

*Universidad de Buenos Aires, Universidade de São Paulo, Universidad
Nacional de Córdoba, Universidad de San Andrés, Universidad de Buenos
Aires and
Universidade de São Paulo*

Efficient automatic protein classification is of central importance in genomic annotation. As an independent way to check the reliability of the classification, we propose a statistical approach to test if two sets of protein domain sequences coming from two families of the Pfam database are significantly different. We model protein sequences as realizations of Variable Length Markov Chains (VLMC) and we use the *context trees* as a signature of each protein family. Our approach is based on a Kolmogorov–Smirnov-type goodness-of-fit test proposed by Balding et al. [Limit theorems for sequences of random trees (2008), DOI: 10.1007/s11749-008-0092-z]. The test statistic is a supremum over the space of trees of a function of the two samples; its computation grows, in principle, exponentially fast with the maximal number of nodes of the potential trees. We show how to transform this problem into a max-flow over a related graph which can be solved using a Ford–Fulkerson algorithm in polynomial time on that number. We apply the test to 10 randomly chosen protein domain families from the seed of Pfam-A database (high quality, manually curated families). The test shows that the distributions of context trees coming from different families are significantly different. We emphasize that this is a novel mathematical approach to validate the automatic clustering of sequences in any context. We also study the performance of the test via simulations on Galton–Watson related processes.

Received May 2008; revised October 2008.

[1]Supported in part by FAPESP, CNPq, Instituto do Milênio for the Global Advance of Mathematics, CAPES-Secyt agreement and PROSUL.

[2]Supported in part by PICT 2005-31659 and PID Secyt 69/08.

[3]Supported by a FAPESP fellowship, proc. 06/56980-0.

**1. Introduction.** The primary structure of a protein is represented by a sequence of 20 different symbols called amino acids. Proteins can be composed of one or more functional regions, called domains; the identification of domains that occur within a protein can provide insights into its function. For this reason biologists classify protein domains into families and care about the reliability of the classification [Stein (2001)]. But in a protein domain database not only the quality of the classification is important, the number of proteins encoded by the genomes that are assigned to the families is also important. This is usually referred to as proteome coverage. For this reason, usually in most databases there must be some balance between quality and quantity.

The Pfam database is a large collection of protein domain families [Finn et al. (2006)]. In its last release of July 2007, the Pfam database comprises 9318 annotated families (Pfam-A) as well as a lower quality, unannotated collection (Pfam-B). Each Pfam-A family consists of two parts: a manually curated set of protein domains called *seed* and a set of automatically detected protein domains using a *profile hidden Markov model* (profile HMM), whose parameters are estimated from the seed of the family.

To our knowledge, no independent method to validate the Pfam classification has been proposed, in spite of problems that the uncertainty in the alignment of sequences can lead to [Wong, Suchard and Huelsenbeck (2008)]. We make a step in this direction by presenting a statistical method to test if two samples from protein domains come from two different families. If some families are not significantly different, then the problem of classifying new proteins becomes risky.

We start by modeling protein sequences as Variable Length Markov Chains (VLMC), a model introduced by Rissanen (1983). A VLMC is a discrete time stochastic process with the property that the law of the process at any given time depends on a finite (but not of fixed length) portion of the process at precedent times [Bühlmann and Wyner (1999)]. As usual in the applications of VLMC, we assume that the process is a Markov chain of order at most $L$ (finite memory process). The minimum set of sequences needed to completely specify the distribution of the next symbol in the sequence is known as a *context tree* and it is denoted by $t$. Calling $p$ the conditional transition probabilities associated to the nodes of $t$, the pair $(t, p)$ completely determines the law of the VLMC.

VLMC have been successfully applied to model and classify protein sequences [Bejerano and Yona (2001)]. As in the case of profile HMM in the construction of the Pfam families, the VLMC approach of Bejerano and Yona takes, for each family, a set of already classified protein domains and estimates a VLMC model, that is, a pair $(t, p)$. Then, the estimated VLMC model is used to classify other protein sequences into the family. Instead, we treat the context trees of sequences of a given family as random samples of a distribution associated to the family; this distribution is used as a *signature* of the family [Galves et al. (2004), Leonardi et al. (2007)]. That is, we propose that the context trees of the sequences disregarding the associated probabilities are sufficient to test if two samples of sequences come from different Pfam families. We take two samples of protein sequences of size $n$ and $m$ respectively and for each sequence we construct the estimated context tree using the PST algorithm introduced by Ron, Singer and Tishby (1996) and implemented by Bejerano (2004), obtaining two samples of trees $\mathbf{t} = (t_1, \ldots, t_n)$, $\mathbf{t}' = (t'_1, \ldots, t'_m)$. We assume that the samples are independent and that the trees in each sample are independent and identically distributed with laws $\pi$ and $\pi'$, respectively. We test $H_0 : \pi = \pi'$ against $H_A : \pi \neq \pi'$ using the test proposed in Balding et al. (2008) (in what follows we will denote it by BFFS test). Rejection of the null hypothesis leads us to conclude that the protein families are distinct.

The BFFS test is a Kolmogorov–Smirnov-type goodness-of-fit test. A distance $d$ defined later in (8) is considered in the space of trees $\mathcal{T}$ and the statistic for the two-sample test is given by

$$(1) \qquad W(\mathbf{t}, \mathbf{t}') := \sup_{t \in \mathcal{T}} |\bar{d}(t, \mathbf{t}) - \bar{d}(t, \mathbf{t}')|,$$

where $\bar{d}(t, \mathbf{t}) = \frac{1}{n} \sum_{i=1}^{n} d(t, t_i)$; that is, $W(\mathbf{t}, \mathbf{t}')$ is the supremum over $t$ in the space of trees $\mathcal{T}$ of the difference of the empiric mean distances of $t$ to each of the two samples $\mathbf{t}$ and $\mathbf{t}'$. The null hypothesis is rejected for large values of $W(\mathbf{t}, \mathbf{t}')$. Since the law of $W$ under $H_0$ is not explicitly known, a simulation procedure is performed to find the $p$-values.

The computation of the test statistic $W(\mathbf{t}, \mathbf{t}')$ is a priori difficult; a naive search would involve an exponential complexity of the algorithm on the number of potential nodes. A major point of this paper is to show that the problem can be re-expressed as to find the maximal flow on a graph constructed as a function of the samples. The approach is inspired by the search for the Maximum a Posteriori in Bayesian image reconstruction using the Ising model, as proposed by Greig, Porteous and Seheult (1989) [see also Kolmogorov and Zabih (2004)], but requires the introduction of a penalty to guarantee that the solution is in $\mathcal{T}$. The max-flow problem can be solved in polynomial time on the maximal number of nodes of the tree, using Ford–Fulkerson type algorithms.

Statistical analysis of tree-like data has been performed in several papers. Banks and Constantine (1998) obtain trees by hierarchical clustering of authors of written texts, using search-related features. They assume a parametric model and use a metric in the space of trees to get a center point and a confidence band around it. Computation of the distribution's parameters, center point and spread are feasible when a distance of the same type as in the BFFS approach is used. Wang and Marron (2007) analyze a sample of blood vessels in the human brain, represented by trees. Each node represents a blood vessel, and the edges represent vessel connections. The statistical analysis of this data was based on a Fréchet approach, which in turn is based on a metric. The Fréchet mean of a data set is the point which minimizes the sum of the squared distances to the data points. In the specific application to blood vessels, both the structure of the trees (i.e., connectivity properties) and the attributes of the nodes (such as the locations and orientations of the blood vessels) were considered. This is a major difference with the BFFS approach, where only the structure of the trees enters in the test statistics.

In Section 2 we define trees, describe VLMC and explain how to obtain the context trees from the observed protein domain sequences. In Section 3 we define the distance in the tree space and describe the BFFS test. In Section 4 we develop the algorithm to compute the BFFS test statistics. In Section 5 we describe the one-sample test and discuss possible extensions of the approach. In Section 6.1 we perform pairwise comparisons of samples of trees corresponding to 10 Pfam families. Final remarks are in Section 7 and computing notes in Section 8. At the end of Section 3 and then in Appendix A.1, we discuss problems related to the power of the tests. Appendix A.2 contains the proofs of selected results. In Appendix A.3 we perform the test on several samples of Galton–Watson related trees obtained with Monte Carlo simulation.

**2. Protein related random trees.** A protein sequence can be modeled as a realization of a discrete time stochastic process having as state space the set $A$ of 20 amino acids. This is the basic idea in the modeling of protein domains by HMMs or VLMCs. In this section we introduce the basic concepts behind VLMC and show how the context tree associated to a protein sequence can be estimated using the Probabilistic Suffix Tree (PST) algorithm proposed by Ron, Singer and Tishby (1996) and implemented by Bejerano (2004).

Let $A$ be a finite alphabet and $V = \bigcup_{\ell=0}^{\infty} A^\ell$ the set of sequences of symbols in $A$. Denote $a_\ell^j$ the sequence $a_\ell a_{\ell+1} \cdots a_j$. Given a sequence $a_1^j$, any sequence $a_\ell^j$ with $1 < \ell \leq j$ is called a *suffix* of $a_1^j$. Let $\mathcal{T} := \{t \subset V : a_1^j \in t \text{ implies } a_2^j \in t\}$ be the space of rooted trees with nodes in $V$; the empty sequence is the root of the tree and it is called $\lambda$. The edges of $t$ are $\{(a_1^j, a_2^j) : a_1^j \in t\}$. A node of an edge is a suffix of the other node and the difference in length of the two nodes is one. Hence, the tree $t$ is identified with its set of nodes.

Let $X = (X_n)_{n \in \mathbb{Z}}$ be a stationary stochastic process taking values in $A$. Define

$$p(a|a_{-j}^{-1}) := P[X_0 = a | X_{-j}^{-1} = a_{-j}^{-1}].$$

A finite sequences $a_{-k}^{-1} \in V$ is sufficient to determine the law of the next symbol if

(2) $\quad p(\cdot|a_{-j}^{-k-1}a_{-k}^{-1}) = p(\cdot|a_{-k}^{-1}) \qquad$ for all $k < j$ and all $a_{-j}^{-k-1} \in A^{j-k}$,

where $a_{-j}^{-k-1}a_{-k}^{-1}$ denotes the concatenation of the sequences $a_{-j}^{-k-1}$ and $a_{-k}^{-1}$. We assume that the process is a Markov chain of order $L$, that is, that (2) holds for $k = L$ and all $a_{-L}^{-1} \in A^L$. A finite sequence $a_{-k}^{-1} \in A^k$ is called a *context* if it satisfies (2) and

(3) $$p(\cdot|a_{-k}^{-1}) \neq p(\cdot|a_{-k+1}^{-1}).$$

We say that $X$ is a VLMC if there are contexts of length less than $L$, that is, if there exists a $k < L$ and a sequence $a_{-k}^{-1}$ satisfying (2) and (3). The set $t$ of contexts and all their suffixes is called the *context tree*; each node in the context tree is labeled by a finite sequence over $A$. In this finite memory case, a VLMC is simply a parsimonious representation of a Markov chain of order $L$ that, in its strict sense, would have $(|A|-1)|A|^L$ parameters (a probability distribution associated to each sequence of fixed length $L$).

Under the assumption of bounded memory, the pair composed by the context tree and the set of transition probability distributions associated to the nodes of $t$ completely specify the law of the stationary process $X$.

Figure 1 summarizes an example of a context tree and transition probabilities for a stationary process $X$ over the alphabet $A = \{1, 2\}$. For $0 < \alpha < 1$ the transition probabilities are given by

(4) $\qquad p(1|a_{-\infty}^{-1}) = \begin{cases} \alpha, & \text{if } a_{-3}^{-1} = 111 \text{ or } a_{-3}^{-1} = 122, \\ 1 - \alpha, & \text{if } a_{-3}^{-1} = 211 \text{ or } a_{-3}^{-1} = 222, \\ 0.5, & \text{otherwise;} \end{cases}$

see Figure 1(a). If $\alpha \neq 0.5$, the set of contexts is $\{111, 122, 211, 222\}$ and the context tree is $t = \{\lambda, 1, 2, 11, 22, 111, 211, 122, 222\}$; see Figure 1(b). If $\alpha = 0.5$, the context tree is just $\lambda$, as the chain is a sequence of i.i.d. (0.5) random variables.

There are several approaches to estimate the context tree and transition probabilities of VLMCs from a finite realization of $X$. We mention the context algorithm proposed by Rissanen (1983); see also Bühlmann and Wyner (1999) and Galves et al. (2008). Recently, Csiszár and Talata (2006) proposed the use of the Bayesian Information Criterion (BIC) and the Minimum Description Length Principle (MDL). These algorithms provide consistent estimates of the parameters. Our work utilizes the PST *algorithm* and so we provide a brief review here.

Suppose $x_1, \ldots, x_l$ is a sample of a VLMC over $A$ specified by the pair $(t, p)$ (in our setting $x_1, \ldots, x_l$ represents a protein over the alphabet of 20

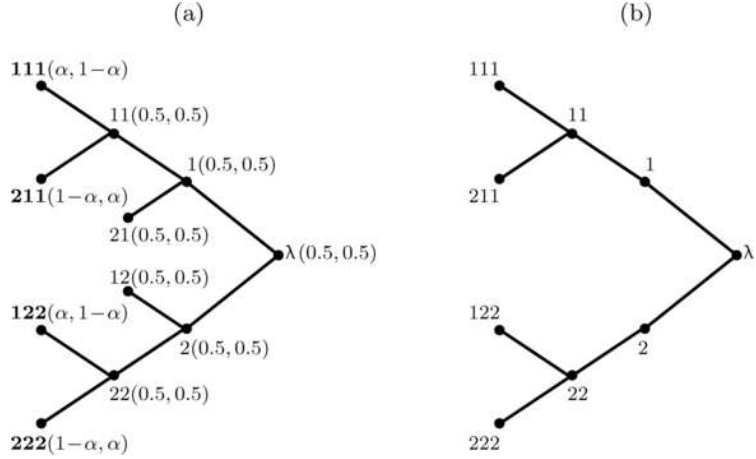

FIG. 1. *An example of stationary conditional probability distributions p over the alphabet $A = \{1, 2\}$ (a) and the corresponding context tree t (b). p and t completely specify a stationary VLMC process $X = (X_n)_{n \in \mathbb{Z}}$. We assume $0 < \alpha < 1$, $\alpha \neq 0.5$. Each node in the tree (a) is labeled by a sequence over the alphabet A and has an associated probability distribution over A (see text for more details). The contexts of the process are the sequences in bold face $\{111, 211, 122, 222\}$ (note that 21 and 12 are not contexts in our definition). The context tree in this case is (b), representing the set $\{\lambda, 1, 2, 11, 22, 111, 211, 122, 222\}$.*

amino acids). For any sequence $a_1^j \in A^j$ define the counters

$$N(a_1^j) = \sum_{i=0}^{l-j} \mathbf{1}\{x_{i+1}^{i+j} = a_1^j\}, \tag{5}$$

where the function **1** takes value 1 if $x_{i+1}^{i+j} = a_1^j$ and 0 otherwise. For any sequence $a_1^j \in A^j$ such that $N(a_1^j) \geq 1$ and any symbol $a \in A$, we define the empirical transition probabilities $\hat{p}(a|a_1^j)$ as

$$\hat{p}(a|a_1^j) = \frac{N(a_1^j a)}{\sum_{b \in A} N(a_1^j b)}. \tag{6}$$

To estimate the context tree associated to the sequence, two parameters are fixed: $L$, the maximal depth of the estimated tree $\hat{t}$ and $r > 1$, a threshold value. The PST algorithm defines the *context tree estimator* $\hat{t}$ as the tree containing all the sequences $a_1^j$, with $j \leq L$, $N(a_1^j) \geq 1$, such that there exists a symbol $a \in A$ satisfying

$$|\log \hat{p}(a|a_2^j) - \log \hat{p}(a|a_1^j)| \geq \log r. \tag{7}$$

That is, the node $a_1^j$ is a node of $\hat{t}$ if the conditional probabilities $\hat{p}(\cdot|a_1^j)$ and $\hat{p}(\cdot|a_2^j)$ are sufficiently far in the sense of (7); this is the empirical version of (2)–(3). To guarantee that $\hat{t}$ is a tree, include also all suffixes of included nodes; that is, $a_1^j \in \hat{t}$ implies $a_2^j \in \hat{t}$.

The PST algorithm uses other parameters to smooth the estimated transition probabilities given by (6). This smoothing is useful to avoid null estimated probabilities that can damage the prediction step in the classification of new sequences. Since our interest is to estimate only the context tree, it is sufficient to consider the parameters $L$ and $r$. We refer the interested reader to Ron, Singer and Tishby (1996), Bejerano (2003) and Bejerano (2004) for a full explanation of the PST algorithm, its implementation and some basic examples.

**3. The tree distance and the two-sample test.** For the two sample problem, the BFFS test is based on the supremum (over the space of trees) of the difference between the empirical mean distance function of each sample to a given tree. More formally, for a node $v \in V$ and a tree $t$ in $\mathcal{T}$ denote $t(v) = \mathbf{1}\{v \text{ is a node of } t\}$. Let $\phi : V \to \mathbb{R}^+$ be a nonnegative function and consider the distance in $\mathcal{T}$ defined by

$$d(t, t') := \sum_{v \in V} \phi(v)(t(v) - t'(v))^2. \tag{8}$$

Let $T$ be a random tree on $\mathcal{T}$ with law $\pi$ and $\mathbf{t} = (t_1, \ldots, t_n)$ a random sample of $T$ (independent random trees with the same law as $T$). Define the empiric expected distance of a tree $t$ to the sample by

$$\bar{d}(t, \mathbf{t}) := \frac{1}{n} \sum_{i=1}^{n} d(t_i, t). \tag{9}$$

Consider two samples $\mathbf{t}$ and $\mathbf{t}'$ of random trees $T$ and $T'$, with laws $\pi$ and $\pi'$, with sample sizes $n$ and $m$ respectively. The two-sample problem is to test

$$H_0 : \pi = \pi', \qquad H_A : \pi \neq \pi'. \tag{10}$$

BFFS show that, under $H_0$, the process

$$\left( \sqrt{\frac{nm}{n+m}} (\bar{d}(t, \mathbf{t}) - \bar{d}(t, \mathbf{t}')), t \in \mathcal{T} \right) \tag{11}$$

converges weakly as $\min(n, m) \to \infty$ to a Gaussian process and propose the statistic

$$W(\mathbf{t}, \mathbf{t}') := \sup_{t \in \mathcal{T}} |\bar{d}(t, \mathbf{t}) - \bar{d}(t, \mathbf{t}')|. \tag{12}$$

Under $H_0$ for large $n$ and $m$, $\sqrt{\frac{nm}{n+m}} W(\mathbf{t}, \mathbf{t}')$ has approximately the law of the supremum over $t$ of a Gaussian process indexed by $t \in \mathcal{T}$. Determining the quantiles $q_\alpha$ using the asymptotic law, the null hypothesis is rejected at level $\alpha$ when

$$|W(\mathbf{t}, \mathbf{t}')| > q_\alpha.$$

The quantiles are obtained using permutation-based randomization techniques [Manly (2007)]. See Section 6.1 for details.

*About the power of the BFFS test.* Strictly speaking, the null hypothesis of the BFFS test is

$$H_0': \text{law of } (d(t,T), t \in \mathcal{T}) = \text{law of } (d(t,T'), t \in \mathcal{T}).$$

Of course, rejection of $H_0'$ implies rejection of $H_0$, so that we do not need to worry when rejecting. But the test could accept $H_0'$ even when $H_0$ is false. We give in Appendix A.1 examples of different $\pi$s giving rise to the same process $(d(t,T), t \in \mathcal{T})$ and show some sufficient conditions on $\pi$ under which the law of $(d(t,T), t \in \mathcal{T})$ determines $\pi$.

**4. Graph computation of $W(\mathbf{t}, \mathbf{t}')$.** To compute the test statistic $W(\mathbf{t}, \mathbf{t}')$, it is necessary to find the trees attaining the supremum (12). In this section we show that the problem can be reformulated in terms of finding the maximal flow on a graph constructed as a function of the sample.

Denote by $\bar{t}$ the empiric mean of the sample $\mathbf{t}$:

$$\bar{t}(v) := \frac{1}{n} \sum_{i=1}^{n} t_i(v). \tag{13}$$

Since $\bar{t}(v)$ is not always an integer, $\bar{t}$ is not necessarily a tree, but $\bar{t}(av) \le \bar{t}(v)$ if $v$ and $av$ are in $V$. Notice that

$$\bar{d}(t, \mathbf{t}) - \bar{d}(t, \mathbf{t}') = 2\mathcal{L}(t) + \sum_{v \in V} \phi(v)(\bar{t}(v) - \bar{t'}(v)), \tag{14}$$

where

$$\mathcal{L}(t) = \sum_{v \in V} \phi(v)(\bar{t'}(v) - \bar{t}(v))t(v). \tag{15}$$

Since the last term in (14) does not depend on $t$, to maximize $|\bar{d}(t, \mathbf{t}) - \bar{d}(t, \mathbf{t}')|$ on $\mathcal{T}$ is equivalent to minimize $\mathcal{L}(t)$ and $-\mathcal{L}(t)$ on $\mathcal{T}$.

Define the space of *configurations*

$$\mathcal{Y} := \{0,1\}^V.$$

This set can be identified with the set of subsets of $V$, so that $\mathcal{T} \subset \mathcal{Y}$. In order to penalize configurations in $\mathcal{Y}$ that are not trees, we define the quadratic function $\mathcal{P} \colon \mathcal{Y} \to \mathbb{N}$ which counts the number of orphan nodes in a configuration:

$$(16) \qquad \mathcal{P}(y) = \sum_{\{v, av\} \subset V} y(av)(1 - y(v)), \qquad y \in \mathcal{Y},$$

where for a node $v = a_1 \cdots a_j$, $av = a a_1 \cdots a_j$ is a son of $v$. It is clear that $\mathcal{P}(y) \geq 0$ and $\mathcal{P}(y) = 0$ if and only if $y \in \mathcal{T}$.

The following proposition shows that to maximize $|\bar{d}(t, \mathbf{t}) - \bar{d}(t, \mathbf{t}')|$ on $\mathcal{T}$ is equivalent to minimize $\mathcal{L}(y) + \beta \mathcal{P}(y)$ and $-\mathcal{L}(y) + \beta \mathcal{P}(y)$ on $\mathcal{Y}$ for $\beta$ sufficiently large.

PROPOSITION 4.1. *Let $\beta > \sum_{v \in V} \phi(v)$. Then*

$$(17) \quad \arg\max_{t \in \mathcal{T}} |\bar{d}(t, \mathbf{t}) - \bar{d}(t, \mathbf{t}')| \\ \subset \arg\min_{y \in \mathcal{Y}}(\mathcal{L}(y) + \beta \mathcal{P}(y)) \cup \arg\min_{y \in \mathcal{Y}}(-\mathcal{L}(y) + \beta \mathcal{P}(y)).$$

The proof of Proposition 4.1 is given in Appendix A.2. The proposition reduces the minimization problem on $\mathcal{T}$ to the task of minimizing the *Hamiltonians* $\mathcal{L} + \beta \mathcal{P}$ and $-\mathcal{L} + \beta \mathcal{P}$ on $\mathcal{Y}$.

The Hamiltonian $\mathcal{L} + \beta \mathcal{P}$ is represented by the oriented graph $(\widetilde{V}, \widetilde{E})$ given by

$$(18) \qquad \widetilde{V} := V \cup \{s\} \cup \{b\}, \qquad \widetilde{E} := \{(v, w) : v, w \in \widetilde{V}\},$$

where $s$ (source) and $b$ (sink) are two extra nodes. The graph has the following capacities associated to the (oriented) edges:

$$(19) \quad c(v, w) := \begin{cases} (\phi(v)(\overline{t}(w) - \overline{t'}(w)))^+, & \text{if } v = s \text{ and } w \in V, \\ (\phi(v)(\overline{t}(v) - \overline{t'}(v)))^-, & \text{if } v \in V \text{ and } w = b, \\ \beta, & \text{if } v \in V, w = av \in V, a \in A, \\ 0, & \text{otherwise}, \end{cases}$$

where $x^+ = \max\{x, 0\}$, $x^- = \max\{-x, 0\}$. That is, the edges linking a node of $V$ to its sons have capacity $\beta$, the edge linking a node of $V$ to the sink, and the edge linking the source to a node of $V$ have capacity $\phi(v)(\overline{t}(v) - \overline{t'}(v)))^\pm$ according to the sign of $(\overline{t}(v) - \overline{t'}(v))$; the other edges have zero capacity.

A configuration $y \in \mathcal{Y}$ defines a *cut* of the graph

$$(20) \quad C(y) := \{(v, w) \in \widetilde{E} : c(v, w) > 0, v \in (V \setminus y) \cup \{s\}, w \in y \cup \{b\}\},$$

whose capacity $c(y)$ is

$$(21) \qquad c(y) := \sum_{(v, w) \in C(y)} c(v, w).$$

The next result has been proven by Kolmogorov and Zabih (2004), as a generalization of an approach of Greig, Porteous and Scheult (1989). We give some details of the proof in Appendix A.2.

PROPOSITION 4.2. *It holds that $c(y) = k + \mathcal{L}(y) + \beta \mathcal{P}(y)$ for all $y \in \mathcal{Y}$, where $k$ does not depend on $y$.*

Proposition 4.2 shows that to minimize $\mathcal{L}(y) + \beta \mathcal{P}(y)$ it is sufficient to find a minimum cut in its associated graph. This problem can be solved by means of the Ford–Fulkerson type of algorithm as proposed by Greig, Porteous and Scheult (1989). We use the variant and implementation of Kolmogorov and Zabih (2004).

The idea behind the Ford–Fulkerson algorithm is the following. Suppose that liquid is flowing from source to sink in the graph with nodes $\widetilde{V}$ and pipes $\widetilde{E}$ with capacities $c(\cdot, \cdot)$. Take a piece of chalk and draw an arrow in the direction of the flow over the pipes with positive capacity that are not totally filled. Draw an arrow in the direction opposite to the flow over the pipes carrying some liquid. Of course there may be pipes with arrows in both directions! Now try to walk from the source to the sink, always following your arrows. When you arrive at a dead end, return to the source. If you never arrive to the sink, the flow is maximal and the nodes $y$ that you have *not* visited define a cut $C(y)$ with minimal capacity. If you arrive to the sink, you can increase the total flow by $\epsilon$ by increasing it by $\epsilon$ in the pipes that you have walked forward from the source, and decreasing it by the same amount in the pipes that you walked backward. It turns out that the number of operations necessary to find the minimal cut is polynomial in the number of nodes.

## 5. Generalizations.

*The one-sample test.* Given a sample $\mathbf{t}$ of a random tree with law $\pi$, the one-sample test is

$$H_0 : \pi = \pi', \qquad H_A : \pi \neq \pi' \tag{22}$$

for a given probability $\pi'$ on $\mathcal{T}$. The BFFS statistic in this case is given by

$$W(\mathbf{t}) := \sup_{t \in \mathcal{T}} |\bar{d}(t, \mathbf{t}) - \pi d(t)|, \tag{23}$$

where

$$\pi d(t) := \sum_{t' \in \mathcal{T}} \pi(t') d(t', t) \tag{24}$$

is the expected distance between $t$ and a random tree with law $\pi$. The graph computation of $W(\mathbf{t})$ is done exactly as in Section 4, but in the definition (15) of $\mathcal{L}$ the mean occupation value $\overline{t'}(v)$ is substituted by

$$\mu_{\pi'}(v) := \sum_{t} \pi'(t) t(v), \tag{25}$$

the average occupation number of node $v$ under $\pi'$.

*More general trees.* The BFFS test works for finite or infinite rooted trees contained in a *full tree* $V$. $V$ must satisfy that the number of children per parent is uniformly bounded by $m$ (say). This condition is necessary for the Proposition 4.1, which transforms the problem of minimizing the difference of the distances on the minimization of a Hamiltonian. The nodes in $V$ can be coded with finite sequences of letters of the alphabet $A = \{1, \ldots, m\}$ in such a way each node is coded with the sequence coding of his father plus a letter of $A$. In our case the letter is added at the beginning of the sequence so that the sequence corresponding to a parent is a suffix of those corresponding to its children. Alternatively, the letter can be added at the end; in this case the parents sequences are prefixes of the children. Any other labeling would work in the same way, as only the structure of the tree (and not the labeling of the nodes) is relevant in the construction of the test. The structure of the tree is the set of nodes and the set of edges; with our coding the set of edges is deduced from the node coding: $E = \{(a_1^j, a_2^j)\}$, otherwise $E$ must be explicitly defined.

If $V$ is infinite, $V$ is truncated to nodes with at most $L$ symbols and calling $\pi_L$ the law of the truncated tree, the null hypothesis is $H_0 : \pi_L = \pi'_L$.

## 6. Numerical results.

6.1. *Testing protein related populations of trees.* In this subsection we present some results obtained by applying the two-sample test over protein domain families of the Pfam-A database. As mentioned in the Introduction, our framework is the following:

- Each family $\mathcal{F}$ of protein domains induce a (different, hopefully) probability distribution $\pi$ on the space of trees $\mathcal{T}$.
- Given two families $\mathcal{F}$ and $\mathcal{F}'$, we consider their associated signatures, that is, the probability laws $\pi$ and $\pi'$ on the space $\mathcal{T}$.
- For each family $\mathcal{F}_j$ we take a sample of protein sequences of size $n_j$, and for each sequence in the sample we construct the PST context tree estimator, as described in Section 2. We obtain a sample of size $n_j$ of i.i.d. random elements on $\mathcal{T}$ with distribution $\pi_j$.
- Finally, for each pair of families $\mathcal{F}_j, \mathcal{F}_{j'}$ we test if both distributions $\pi_j$ and $\pi_{j'}$ are the same.

To test the approach, we randomly choose families $\mathcal{F}_1, \ldots, \mathcal{F}_{10}$ whose names start with letter A and such that their average lengths are larger than 150 amino acids (this last condition was to ensure some precision in the context tree estimation step). In order to guarantee the quality of the samples, we only choose sequences in the seed of each family. The chosen families are ABC-2membrane, ABC-membrane, Amidase, Amidino-trans, AMME-CR1, AOX, ArgK, ASC, Asp-Arg-Hidrox and Asp-Al-Ex.

We randomly select $n_j = 50$ sequences from each family $\mathcal{F}_j$ and compute the associated PST context tree estimator of each sequence using the PST algorithm with parameters $r = 1.05$ [as Bejerano (2001)] and $L = 4$. In this way we obtain a sample of 50 trees per family.

We consider the distance function $\phi(v) = \theta^{\text{gen}(v)}$, where the function $\text{gen}(v)$ is defined as the length of the sequence labeling node $v$. That is, if node $v$ corresponds to sequence $a_1^k$, then $\text{gen}(v) = k$. For each $\theta \in \{0.001, 0.01, 0.35\}$ we run the BFFS test for each pair of families using the corresponding samples of trees. We also run the tests under the null hypothesis collecting two independent samples from the same family. For each pair of samples of trees we estimate the $(1-\alpha)$-quantile under the null hypothesis using Monte Carlo randomization [Manly (2007)], that is, we permute the pooled sample a thousand times and compute the test statistic for each of the replicates using half of the permuted pooled sample for each population. The estimated quantile is therefore the empirical $(1-\alpha)$-quantile for the vector of size 1000 built up in this way.

All 45 tests are rejected at level $\alpha = 0.001$ for the three values of $\theta$. Despite the conservative level used in the tests, the hypothesis of equal distribution is rejected in all cases when the samples came from different Pfam-A families, confirming the discriminative power of the context trees associated to the sequences. In the case of the same family, but with independent samples of trees, for $\theta = 0.35$ we observe $p$-values ranging from 0.15 to 0.87, compatible with the uniform distribution (the law of the $p$-values under the null hypothesis). Similar results were obtained with the other two values of $\theta$.

6.2. *Simulation results.* We also challenge our method in a small Monte Carlo simulation for Galton–Watson processes, for three different models, parameters and sample sizes. The results are reported in Appendix A.3.

**7. Final remarks.** We perform the BFFS method to test if two samples of context trees come from different distributions, and we propose a feasible way to compute its statistic, allowing the treatment of reasonably big trees. The test rejects the null hypothesis in the case of high quality, manually curated Pfam families, and it does not reject on random subsets of the same family. This supports the use of the test as a method to distinguish different groups of protein domains when a specific task, as, for example, sequence annotation, does not give conclusive results.

Our results strongly indicate that the context trees associated to protein domain sequences are sufficient to discriminate between different families in the Pfam database. In this sense we have benefited from ideas coming from the analysis of sequences related to linguistics. Galves and collaborators have proposed with success the use of context trees to discriminate languages from codified written text [Galves et al. (2004)]. More recently, they have used similar ideas in a preliminary work to study the phylogeny of protein sequences [Leonardi et al. (2007)].

We emphasize that the test is not restricted to the analysis of samples of context trees. Any space of trees satisfying the assumptions of Section 5 will be suitable for using our approach. On the other hand, for particular distributions like Galton–Watson processes, more simple tailor-made tests can be developed. Our simulations show that the BFFS test is able to distinguish between distributions determined by the node-marginal distributions, which is a large family of distributions for applications. This class includes tree laws with a Markovian hypothesis, as shown in Proposition A.2. We discuss this item in detail in Appendix A.1.

**8. Computing notes.** The code to compute the test statistic is available from Jorge R. Busch (jbusch@fi.uba.ar) upon request. Calculation reported here used Scilab INRIA (http://www.scilab.org/) and C++ code from Bejerano (2004) and Kolmogorov and Zabih (2004).

The computational burden for our algorithm allows us to work with trees with up to $3^{20}$ nodes. Each $p$-value involves 1001 test statistic calculations, with a sample size $n = 50$ for each family, taking at most 15 minutes to be complete, with a Pentium Core 2 duo with 2Gb of RAM memory.

## APPENDIX

### A.1. Mean distances and Markovian hypotheses.

*Do the mean distances determine a measure?* Recall the mean distance $\pi d$ is defined in (24). Our test is universally consistent within the class of distributions for which $(\pi d(t), t \in \mathcal{T})$ determine the probability $\pi$. That is, whatever is the law of the families, the test will asymptotically detect the difference. We show in Lemma A.1 that $\pi d$ determines the law of the marginals $(T(v), v \in V)$. Proposition A.2 says that, under Markov type hypotheses, the marginal distributions determine the measure.

For a random tree $T$ with law $\pi$ recall $\mu_\pi(v)$ is the mean occupancy node $v$ defined in (25) and define $\sigma_\pi^2(v)$, the variance of $T(v)$, by

(A.1) $$\sigma_\pi^2(v) = \mu_\pi(v)(1 - \mu_\pi(v)).$$

LEMMA A.1. *Let $\pi$ and $\pi'$ be measures on $\mathcal{T}$. Then, $\pi d = \pi' d$ if and only if $\mu_\pi = \mu_{\pi'}$.*

The proof is given later in this section.

Different measures may have the same mean distances. For instance, consider $A = \{1, 2\}$ and $\pi$, $\pi'$ defined by

$$\pi(\varnothing) = \tfrac{1}{2}; \qquad \pi(\{\lambda\}) = \pi(\{\lambda, 1\}) = \pi(\{\lambda, 2\}) = \pi(\{\lambda, 1, 2\}) = \tfrac{1}{8},$$

$$\pi'(\varnothing) = \tfrac{1}{2}; \qquad \pi'(\{\lambda\}) = \pi'(\{\lambda, 1, 2\}) = \tfrac{1}{4}.$$

Then, $\mu_\pi(v) = \mu_{\pi'}(v)$ for all $v \in V$.

Lemma A.1 implies that in general the functions $\pi d$ and $\pi' d$ do not help to solve the discrimination problem. But in some cases these functions do discriminate. To show that, we need some extra notation. For a set $I$ of nodes denote $T_I$ the restriction of $T$ to $I$ and $T_I = 1$ means that $T(v) = 1$ for all $v \in I$, while $T_I = 0$ means that $T(v) = 0$ for all $v \in I$.

Let $v$ be a node, and $a, b \in A$. We shall call $v$ *father* of $av$, $av$ *son* of $v$, and $av$ *brother* of $bv$. Let $f : V \setminus \{\lambda\} \to V$ be a function such that, for each $v \neq \lambda$, $f(v)$ is father or brother of $v$, and $f^n(v) = \lambda$ for some $n = n(v) \in \mathbb{N}$. Notice that, in this case, $f^{-1}(v)$ is empty or formed up by brothers and sons of $v$. We call such a function a *tree-shift*. Consider, for instance, the function $f$ that assigns to a node its father.

Let $T$ be a random tree with law $\pi$. We say that $T$ satisfies a *Markov hypothesis* if there exists a tree-shift $f$ such that if $v = f(w)$, then

(a) $0 < \mu_\pi(w) < \mu_\pi(v) < 1$,
(b) $\mathbb{P}(T(w) = 1 | T(v) = 0) = 0$,
(c) let $I, J \subset V$ be such that if $v \in I$, $w \notin I \cup J$ and $\mathbb{P}(T_I = 1, T_J = 0) > 0$,

then

(A.2) $\qquad \mathbb{P}(T(w) = 1 | T_I = 1, T_J = 0) = \mathbb{P}(T(w) = 1 | T(v) = 1)$.

PROPOSITION A.2. *Under the Markov hypotheses, the marginals $(\mu_\pi(v), v \in V)$ determine the probabilities $(\pi(t), t \in \mathcal{T})$.*

The proof is given later in this section.

*Examples of measures satisfying the Markov hypotheses.* The alphabet for the following examples is $A = \{1, \ldots, m\}$.

1. Let $f$ be defined by $f(av) = v$, $a \in A$. Let $k(v)$ be such that $f^{k(v)-1}(v) = 1$, for $v \in V$. If $\mu_\pi(v) = p^{k(v)}$ ($0 < p < 1$), we obtain a tree with $\pi(\{\lambda\}) = p$, and when $T(v) = 1$, $T(av)$ is Bernoulli with parameter $p$, for $a \in A$. We call such tree distributions *pseudo Galton–Watson processes*.

2. Let $f$ be defined by $f(1v) = v$, and $f((a+1)v) = av$ for $1 \leq a < m$. That is, $v = f(w)$ if $w$ is the eldest brother and $v$ is the father of $w$, or if $v$ is the nearest older brother of $v$. Let now $p_0, \ldots, p_m$ be given probabilities, with $p_0 > 0$ and $p_0 + \cdots + p_m = 1$. If

$$\mu_\pi(1v) = (p_1 + \cdots + p_m)\mu_\pi(v),$$
$$\mu_\pi((a+1)v) = \frac{p_{a+1} + \cdots + p_m}{p_a + \cdots + p_m}\mu_\pi(av) \qquad (1 \leq a \leq m-1),$$

we obtain the classical Galton–Watson process, with parameter probabilities $p_0, \ldots, p_m$.

If $\pi$ is a distribution on $\mathcal{T}$, and $T$ is a random tree with distribution $\pi$, then by a simple computation,

$$\pi d(t) = \sum_{v \in V} \phi(v)(\mu_\pi(v) - t(v))^2 + \sum_{v \in V} \phi(v)\sigma_\pi^2(v) \tag{A.3}$$

$$= \sum_{v \in V} \phi(v)\mu_\pi(v)(1 - 2t(v)) + \sum_{v \in V} \phi(v)t(v). \tag{A.4}$$

PROOF OF LEMMA A.1. Notice first that from (A.4) it follows that

$$\pi d(t) - \pi' d(t) = \sum_{v \in V} \phi(v)(\mu_\pi(v) - \mu_{\pi'}(v))(1 - 2t(v)), \tag{A.5}$$

which implies that if $\mu_\pi(v) = \mu_{\pi'}(v)$ for all $v \in V$, then $\pi d(t) = \pi' d(t)$. This proves sufficiency.

To prove necessity, we proceed by induction. When $\pi d(t) = \pi' d(t)$ for all $t \in \mathcal{T}$, from (A.5) we obtain

$$\begin{aligned}
0 &= \sum_{v \in V} \phi(v)(\mu_\pi(v) - \mu_{\pi'}(v))(1 - 2t(v)) \\
&= -\sum_{v \in t} \phi(v)(\mu_\pi(v) - \mu_{\pi'}(v)) + \sum_{v \notin t} \phi(v)(\mu_\pi(v) - \mu_{\pi'}(v))
\end{aligned} \tag{A.6}$$

for all $t \in \mathcal{T}$. Letting $t = \varnothing$, the empty tree, and $t = \{\lambda\}$ in (A.6), we obtain

$$0 = \sum_{v \in V} \phi(v)(\mu_\pi(v) - \mu_{\pi'}(v)), \tag{A.7}$$

$$0 = -\phi(\lambda)(\mu_\pi(\lambda) - \mu_{\pi'}(\lambda)) + \sum_{v \neq \lambda} \phi(v)(\mu_\pi(v) - \mu_{\pi'}(v)). \tag{A.8}$$

Substracting (A.7) and (A.8) and using that $\phi(v) > 0$, we get $\mu_\pi(\lambda) = \mu_{\pi'}(\lambda)$.

*Inductive step.* Let $t \in \mathcal{T}$, and $h \in V \setminus \{t\}$ such that $t \cup \{h\} \in \mathcal{T}$. We show that if $\mu_\pi(v) = \mu_{\pi'}(v)$ for all $v \in t$, then $\mu_\pi(h) = \mu_{\pi'}(h)$. First, we obtain from (A.6)

$$0 = \sum_{v \notin t} \phi(v)(\mu_\pi(v) - \mu_{\pi'}(v)), \tag{A.9}$$

$$0 = -\phi(h)(\mu_\pi(h) - \mu_{\pi'}(h)) + \sum_{v \notin t \cup \{h\}} \phi(v)(\mu_\pi(v) - \mu_{\pi'}(v)), \tag{A.10}$$

and it follows that $\mu_\pi(h) = \mu_{\pi'}(h)$.  $\square$

It is easy to prove the following lemma

LEMMA A.3. *If $T$ is a random tree satisfying the Markov hypotheses with tree-shift $f$, then, given $T(v) = 1$, the variables $T(w): w \in f^{-1}(v)$ are independent. Furthermore, if $v = f(w)$,*

$$\mathbb{P}(T(w) = 1 | T(v) = 1) = \frac{\mu_\pi(w)}{\mu_\pi(v)}. \tag{A.11}$$

PROOF OF PROPOSITION A.2. First, notice that

$$\pi(\varnothing) = 1 - \mu_\pi(\lambda). \tag{A.12}$$

From Lemma A.3,

$$\pi(\{\lambda\}) = \mu_\pi(\lambda) \prod_{h \in f^{-1}(\lambda)} \left(1 - \frac{\mu_\pi(h)}{\mu_\pi(\lambda)}\right). \tag{A.13}$$

Let $t \in \mathcal{T} \setminus \{\varnothing, V\}$ and $h$ be a node such that $h \notin t$ and $v = f(h) \in t$. We shall show that

$$\pi(t \cup \{h\}) = \pi(t) \frac{\mu_\pi(h)}{\mu_\pi(v) - \mu_\pi(h)}. \tag{A.14}$$

First, we have

$$\begin{aligned}
\pi(t \cup \{h\}) &= \mathbb{P}(T_t = 1, T(h) = 1, T_{(t \cup \{h\})^c} = 0) \\
&= \mathbb{P}(T(h) = 1 | T_t = 1, T_{(t \cup \{h\})^c} = 0)\mathbb{P}(T_t = 1, T_{(t \cup \{h\})^c} = 0) \\
&= \mathbb{P}(T(h) = 1 | T(v) = 1)\mathbb{P}(T_t = 1, T_{(t \cup \{h\})^c} = 0).
\end{aligned} \tag{A.15}$$

On the other hand,

$$\begin{aligned}
\pi(t) &= \mathbb{P}(T_t = 1, T(h) = 0, T_{(t \cup \{h\})^c} = 0) \\
&= \mathbb{P}(T(h) = 0 | T(v) = 1)\mathbb{P}(T_t = 1, T_{(t \cup \{h\})^c} = 0) \\
&= (1 - \mathbb{P}(T(h) = 1 | T(v) = 1))\mathbb{P}(T_t = 1, T_{(t \cup \{h\})^c} = 0).
\end{aligned} \tag{A.16}$$

From (A.15) and (A.16) it follows that

$$\pi(t \cup \{h\}) = \pi(t) \frac{\mathbb{P}(T(h) = 1 | T(f(h)) = 1)}{1 - \mathbb{P}(T(h) = 1 | T(f(h)) = 1)}. \tag{A.17}$$

This shows (A.14). Our main statement follows now by induction from (A.13) and (A.14), noticing that any finite tree may be constructed from $\{\lambda\}$ in this way. □

**A.2. Redefining the minimization problem.** In this subsection we prove Propositions 4.1 and 4.2. Call $\Delta(v) = (\overline{t}(v) - \overline{t'}(v))$ so that

$$\mathcal{L}(t) = \sum_{v \in V} \phi(v) \Delta(v) t(v). \tag{A.18}$$

Recall the space of configurations $\mathcal{Y} = \{0,1\}^V$. Trees minimizing $\mathcal{L}$ will also minimize $\mathcal{L} + \beta \mathcal{P}$ for all positive $\beta$: If $t \in \mathcal{T} \cap \arg\min_{y \in \mathcal{Y}} \mathcal{L}(y)$, then $t \in \arg\min_{y \in \mathcal{Y}} (\mathcal{L}(y) + \beta \mathcal{P}(y))$ for all $\beta > 0$. On the other hand, if $\beta$ is big enough, we expect the configurations minimizing $\mathcal{L} + \beta \mathcal{P}$ to be trees. Since on the set of trees the form $\mathcal{P}$ vanishes, the minimizing trees should also minimize $\mathcal{L}$. This is proven in the following lemma.

LEMMA A.4. *If $\beta > \sum_{v \in V} \phi(v)$, then*

$$\arg\min_{t \in \mathcal{T}} \mathcal{L}(t) = \arg\min_{y \in \mathcal{Y}} (\mathcal{L}(y) + \beta \mathcal{P}(y)). \tag{A.19}$$

PROOF. Observe that minimizing configurations in $\mathcal{Y}$ satisfy

$$y' \in \arg\min_{y \in \mathcal{Y}} \mathcal{L}(y) \quad \text{if and only if} \quad y'(v) = \begin{cases} 0, & \text{if } \Delta(v) > 0, \\ 1, & \text{if } \Delta(v) < 0. \end{cases} \tag{A.20}$$

Let $\mathcal{L}_{\min}$ and $\mathcal{L}_{\max}$ be the values of the minimum and maximum of $\mathcal{L}$ over $\mathcal{Y}$, that is,

$$\mathcal{L}_{\min} = \sum_{v : \Delta(v) < 0} \phi(v) \Delta(v), \qquad \mathcal{L}_{\max} = \sum_{v : \Delta(v) > 0} \phi(v) \Delta(v). \tag{A.21}$$

Notice that

$$\mathcal{L}_{\max} - \mathcal{L}_{\min} = \sum_{v \in V} \phi(v) |\Delta(v)| \leq \sum_{v \in V} \phi(v). \tag{A.22}$$

For all $y \in \mathcal{Y}$ it holds

$$\mathcal{L}(y) + \beta \mathcal{P}(y) \geq \mathcal{L}_{\min} + \beta \sum_{v : y(v) = 0} \sum_{a \in A} y(av)(1 - y(v)). \tag{A.23}$$

If $y$ is not a tree, there exists $v \in V$ and $a \in A$ such that $y(v) = 0$ and $y(av) = 1$, hence,

$$\mathcal{L}(y) + \beta \mathcal{P}(y) \geq \mathcal{L}_{\min} + \beta > \mathcal{L}_{\min} + \sum_{v \in V} \phi(v) \geq \mathcal{L}_{\max} \tag{A.24}$$

by (A.22). On the other hand, $\mathcal{L}(t) + \beta \mathcal{P}(t) = \mathcal{L}(t) \leq \mathcal{L}_{\max}$ for any tree $t$. Hence, for these values of $\beta$, if $y$ is not a tree, then $\mathcal{L}(y) + \beta \mathcal{P}(y) > \max_{t \in \mathcal{T}} (\mathcal{L}(t) + \beta \mathcal{P}(t))$ and the result follows. □

PROOF OF PROPOSITION 4.1. It follows from (14), (15) and the above lemma applied to $\mathcal{L}$ and $-\mathcal{L}$. □

PROOF OF PROPOSITION 4.2. The proof follows from a simple algebra. More generally, the Hamiltonian $\mathcal{H}: \mathcal{Y} \to \mathbb{R}$ given by

$$\text{(A.25)} \qquad \mathcal{H}(y) = \sum_{v \in V} H^v(y(v)) + \sum_{v,w \in V} H^{v,w}(y(v), y(w))$$

is said *regular* if the quadratic terms satisfy

$$\text{(A.26)} \qquad H^{v,w}(0,0) + H^{v,w}(1,1) \leq H^{v,w}(0,1) + H^{v,w}(1,0).$$

Theorem 4.1 of Kolmogorov and Zabih (2004) says that regular Hamiltonians are *graph representable*, meaning that it is possible to associate capacities to the graph $(\widetilde{V}, \widetilde{E})$ defined in (18) in such a way that

$$\text{(A.27)} \qquad C(y) = k + \mathcal{H}(y), \qquad y \in \mathcal{Y},$$

where $C(y)$ is the capacity of the cut defined by $y$ [see (21)] and $k$ is a constant.

The Hamiltonian $\mathcal{L} + \beta \mathcal{P}$ is regular because $H^{v,av}(t(v), t(av)) = \beta t(av)(1 - t(v))$ satisfy (A.26). The graph (18) with capacities (19) is the Kolmogorov and Zabih representation of the Hamiltonian $\mathcal{L} + \beta \mathcal{P}$. □

**A.3. Simulation results.** In this subsection we perform the test for the two-sample problem (10) using samples $\mathbf{t}$ and $\mathbf{t}'$ with distribution $\pi$ and $\pi'$ on $\mathcal{T}$, trees with maximal depth $L = 8$. We compute the BFFS statistic $W$ given in (12) using the approach of Section 4. We use the distance (8) with $\phi(v) = \theta^{\text{gen}(v)}$ for various values of $\theta$. Since in this case we know $\pi$ and $\pi'$, we estimate the quantiles directly by Monte Carlo simulation, as follows:

Quantile
1. Generate two samples of size $n$ both from $\frac{1}{2}\pi + \frac{1}{2}\pi'$, a fair mixture of the laws. Label them *sample 1* and *sample 2*. Compute the test statistic $W$ using the samples.
2. Repeat the above procedure a fixed number of times $N$.
3. Order the computed statistics values increasingly and define the quantile $q(1 - \alpha)$ as the statistic in place $(1 - \alpha)N$.
4. Calculate the quantiles for several values of $\alpha$ and sample sizes $n$.

Power
1. Generate sample 1 from $\pi$, sample 2 from $\pi'$ and compute $W$ using them.
2. Compare the obtained value against the quantile, and reject the null hypothesis with level $\alpha$ if $W > q(1 - \alpha)$.
3. Repeat the last two steps a fixed number of times and compute the percentage of rejections for each value of $\alpha$ as a measure of the power of the test.

TABLE 1
*Model 1. Power of the tests with $p = 0.5$ and $p' = 0.6$, 0.7, 0.8, sample size $n = 31$, 51, 125*

| $\alpha/p'$ | $n = 31$ | | | $n = 51$ | | | $n = 125$ | | |
| --- | --- | --- | --- | --- | --- | --- | --- | --- | --- |
| | 0.6 | 0.7 | 0.8 | 0.6 | 0.7 | 0.8 | 0.6 | 0.7 | 0.8 |
| 0.01 | 0.101 | 0.504 | 0.923 | 0.122 | 0.744 | 0.999 | 0.466 | 0.994 | 0.9999 |
| 0.05 | 0.242 | 0.708 | 0.986 | 0.313 | 0.889 | 0.9999 | 0.661 | 0.999 | 0.9999 |
| 0.1 | 0.325 | 0.806 | 0.998 | 0.426 | 0.935 | 0.9999 | 0.743 | 0.9999 | 0.9999 |

MODEL 1: BINOMIAL. Let $\pi$ be the law of a Galton–Watson process with offspring distribution Binomial$(2,p)$ and $\pi'$ is the same with parameter $p'$. We use $p = 0.5$ and $p' = 0.6$, 0.7, 0.8. Table 1 shows the percentage of rejection over 1000 tests of level $\alpha = 0.10$, 0.05, 0.01 for sample sizes $n = 31, 51, 125$. The results show the consistency of the BFFS test for alternatives with any value of $p \neq 0.5$. For this simple model, small sample sizes are enough to get high power.

MODEL 2: MIXTURE OF BINOMIALS. Let $\pi$ be the law of a Galton–Watson process with offspring distribution a mixture of Binomials. Independently at each node with probability $q$ use a Binomial$(2, p_1)$, otherwise a Binomial$(2, p_2)$. For $\pi'$ we use $q'$, $p'_1$ and $p'_2$. But we take $q = q' = 0.5$ in the examples.

MODEL 2.1. Take $p_1 = 0.45$, $p_2 = 0.5$ and $p'_1 = 0.1$, $p'_2 = 0.85$. Figure 2 shows histograms of the test statistic values obtained in 1000 iterations, for sample sizes 31 (left) and 131 (right) respectively. We have plotted the test statistic under null hypothesis on red and the one under the alternative on blue. The distance parameter was $\theta = 0.35$. The expected mean at each node is the same because $p_1 + p_2 = p'_1 + p'_2$; the distributions under the null and alternative hypothesis are close to each other but the BFFS test needs only a moderate sample size to give high power to the test. In Table 2 we evaluate the power of the test for three different values of the test level ($\alpha$, as 0.01, 0.05 and 0.1), in the same way as in Table 1. The left plot on Figure 3 shows 5 curves of 1000 $p$-values each, for sample sizes $N = 31, 51, 71, 101, 131$, and the distance parameter $\theta = 0.35$. The right one shows the same results when the distance parameter is changed to $\theta = 0.49$. Increasing the parameter $\theta$ decreases the power of the test.

MODEL 2.2. The parameters $p_1 = 0.3$, $p_2 = 0.65$, $p'_1 = 0.45$ and $p'_2 = 0.5$ give a different scenario, since the distributions of the populations are very close. Figure 4 shows slow changes in the empirical distribution as the sample size grows. At the left is the histogram of the 1000 values of the test statistic under the null (red plot) and alternative hypothesis (blue plot) for sample size $N = 50$. For the right histogram the sample size is $N = 500$. In Figure 5 we consider larger sample sizes for Model 2.2. We fix $\theta = 0.35$ and plot the $p$-values computed over 1000 replications for sample sizes $N = 50$, 250 and 500 (left plot) and the percentage of rejection as a function of the sample size, each curve computed with a different $\alpha$ level. Besides sample fluctuations, the percentage of rejection increases with sample size, but at a quite slow rate.

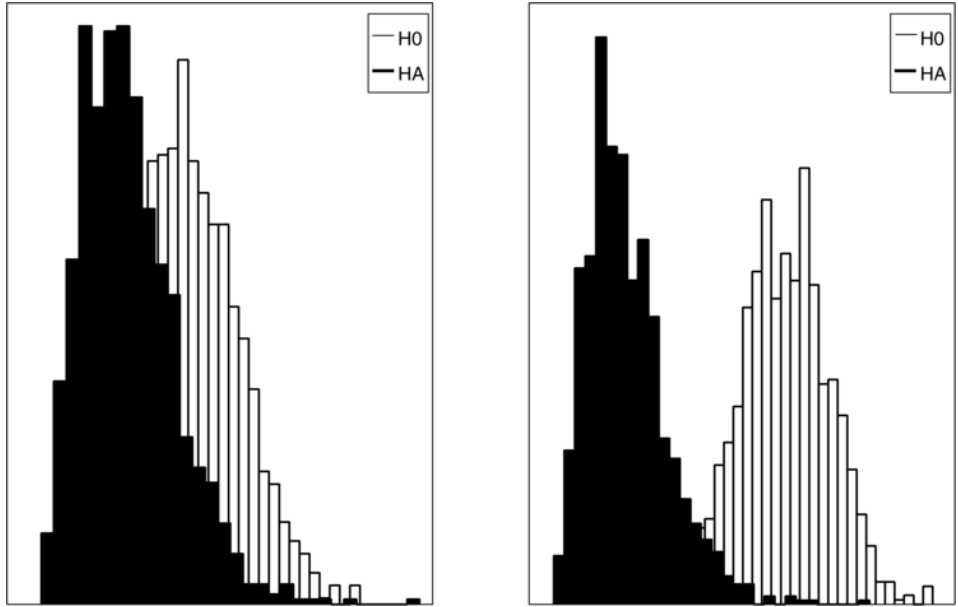

Fig. 2. *Model 2.1. Histogram of the test statistic under null hypothesis (black plot) and alternative (white plot), with sample sizes 31 and 131.*

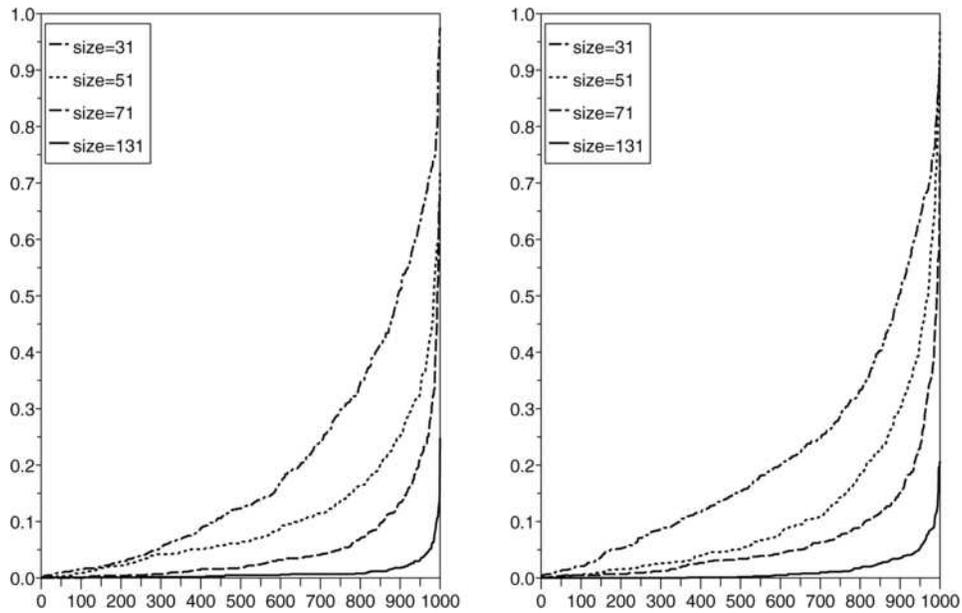

Fig. 3. *Model 2.1. p-values calculated 1000 times; sample sizes 31, 51, 71 and 131. Left plot: parameter $\theta = 0.35$. Right plot: parameter $\theta = 0.49$. Increasing the parameter $\theta$ decreases the power of the test.*

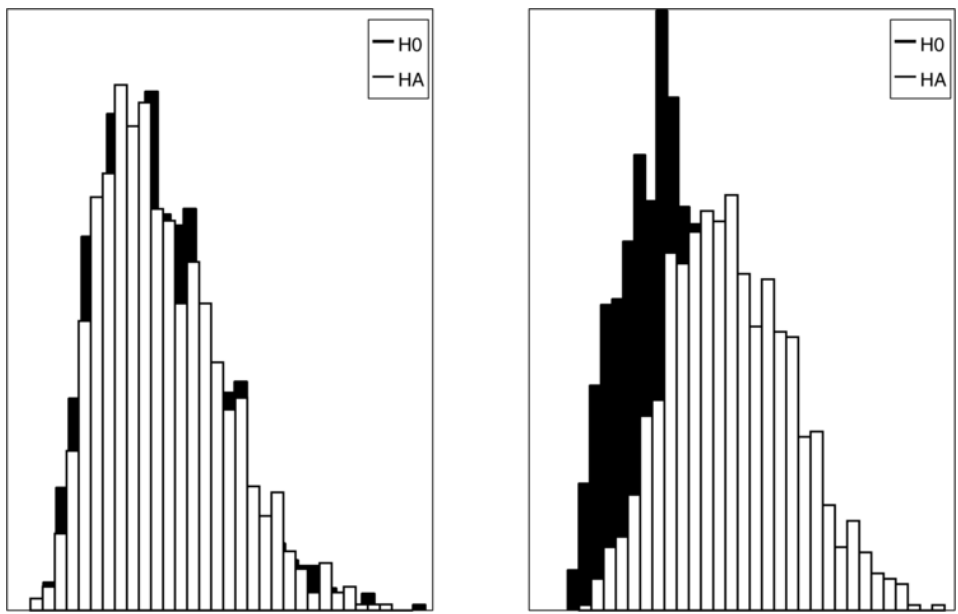

Fig. 4. *Model 2.2. Histogram of the test statistic under null hypothesis (black plot) and alternative (white plot); left plot size = 50 and right plot size = 500. Parameter $\theta = 0.35$.*

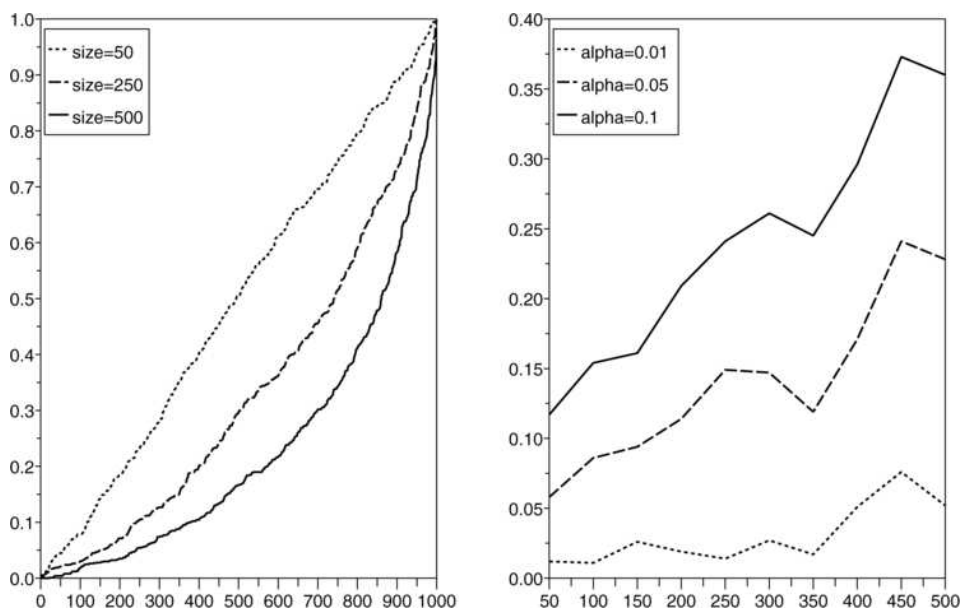

Fig. 5. *Model 2.2. Left plot: p-values computed 1000 times, each curve related to a different sample size 50, 250 and 500. Right plot: percentage of rejection as a function of the sample size, each curve computed with a different $\alpha$ level, 0.01 in dotted line, 0.05 in dashed line and 0.1 in interpolated line. Parameter $\theta = 0.35$.*

TABLE 2
*Model 2.1. Power of the tests of level $\alpha$ for 1000 replications, with $\alpha = 0.05$, 0.01, 0.1. Sample sizes are 31, 51, 71, 101, 131. Parameter $\theta = 0.35$*

| $\alpha$ | 31    | 51    | 71    | 101   | 131   |
|----------|-------|-------|-------|-------|-------|
| 0.01     | 0.045 | 0.108 | 0.325 | 0.594 | 0.819 |
| 0.05     | 0.288 | 0.363 | 0.747 | 0.887 | 0.973 |
| 0.1      | 0.438 | 0.642 | 0.857 | 0.976 | 0.990 |

**Acknowledgments.** We thank Antonio Galves for illuminating discussions about the discriminative power of context trees.

We thank the referees, an associate editor and the editor of the journal for their careful reading of the manuscript and several comments which helped to improve the paper.

## REFERENCES

BALDING, D., FERRARI, P., FRAIMAN, R. and SUED, M. (2008). Limit theorems for sequences of random trees. *Test* (online). DOI: 10.1007/s11749-008-0092-z.

BANKS, D. and CONSTANTINE, G. (1998). Metric models for random graphs. *J. Classification* **15** 199–223. MR1665974

BEJERANO, G. (2003). Automata learning and stochastic modeling for biosequence analysis. Ph.D. thesis, Hebrew Univ.

BEJERANO, G. (2004). Algorithms for variable length Markov chain modeling. *Bioinformatics* **20** 788–789.

BEJERANO, G. and YONA, G. (2001). Variations on probabilistic suffix trees: Statistical modeling and prediction of protein families. *Bioinformatics* **17** 23–43.

BÜHLMANN, P. and WYNER, A. J. (1999). Variable length Markov chains. *Ann. Statist.* **27** 480–513. MR1714720

CSISZÁR, I. and TALATA, Z. (2006). Context tree estimation for not necessarily finite memory processes, via BIC and MDL. *IEEE Trans. Inform. Theory* **52** 1007–1016. MR2238067

FINN, R. D., MISTRY, J., SCHUSTER-BÖCKLER, B., GRIFFITHS-JONES, S., HOLLICH, V., LASSMANN, T., MOXON, S., MARSHALL, M., KHANNA, A., DURBIN, R., EDDY, S. R.,

SONNHAMMER, E. L. L. and BATEMAN, A. (2006). Pfam: Clans, web tools and services. *Nucleic Acids Res.* **34** D247-D51.

GALVES, A., GALVES, C., GARCIA, N. and PEIXOTO, C. (2004). Correlates of rhythm in written texts of Brazilian and European Portuguese. Preprint. Available as [Technical Report 08/08, IMECC/UNICAMP](Technical Report 08/08, IMECC/UNICAMP).

GALVES, A., MAUME-DESCHAMPS, V. and SCHMITT, B. (2008). Exponential inequalities for VMLC empirical trees. *ESAIM Probab. Stat.* **12** 219–229. MR2374639

GREIG, D., PORTEOUS, B. and SEHEULT, A. (1989). Exact maximum a posteriori estimation for binary images. *J. Roy. Statist. Soc. Ser. B* **51** 271–279.

KOLMOGOROV, V. and ZABIH, R. (2004). What energy functions can be minimized via graphs cuts? *IEEE Trans. Pattern Analysis and Machine Intelligence* **26** 147–159.

LEONARDI, F., MATIOLI, S. R., ARMELIN, H. A. and GALVES, A. (2007). Detecting phylogenetic relations out from sparse context trees. Available at [ArXiv:math/0804.4279](ArXiv:math/0804.4279).

MANLY, B. F. J. (2007). *Randomization, Bootstrap and Monte Carlo Methods in Biology*, 3rd ed. Chapman & Hall/CRC, New York.

RISSANEN, J. (1983). A universal data compression system. *IEEE Trans. Inform. Theory* **29** 656–664. MR0730903

RON, D., SINGER, Y. and TISHBY, N. (1996). The power of amnesia: Learning probabilistic automata with variable memory length. *Machine Learning* **25** 117–149.

STEIN, L. (2001). Genome annotation: From sequence to biology. *Nat. Rev. Genet.* **2** 493–505.

WANG, H., and MARRON, J. S. (2007). Object oriented data analysis: Sets of trees. *Ann. Statist.* **35** 1849–1873. MR2363955

WONG, K. M., SUCHARD, M. A. and HUELSENBECK, J. P. (2008). Alignment uncertainty and genomic analysis. *Science* **319** 473–476. MR2381044

J. R. BUSCH
S. P. GRYNBERG
DEPARTAMENTO DE MATEMÁTICAS
FACULTAD DE INGENIERÍA
UNIVERSIDAD DE BUENOS AIRES
ARGENTINA
E-MAIL: [jbusch@fi.uba.ar](jbusch@fi.uba.ar)
   [sebgryn@fi.uba.ar](sebgryn@fi.uba.ar)

P. A. FERRARI
F. LEONARDI
INSTITUTO DE MATEMÁTICA E ESTATÍSTICA
UNIVERSIDADE DE SÃO PAULO
BRAZIL
E-MAIL: [pablo@ime.usp.br](pablo@ime.usp.br)
   [florencia@usp.br](florencia@usp.br)

A. G. FLESIA
CENTRO DE INVESTIGACIÓN Y ESTUDIOS
DE MATEMÁTICA DE CÓRDOBA, CONICET
AND
FAMAF-UNC, ING. MEDINA ALLENDE S/N
CIUDAD UNIVERSITARIA
CP 5000, CÓRDOBA
ARGENTINA
E-MAIL: [flesia@mate.uncor.edu](flesia@mate.uncor.edu)

R. FRAIMAN
DEPARTAMENTO DE MATEMÁTICA Y CIENCIAS
UNIVERSIDAD DE SAN ANDRÉS
BUENOS AIRES
ARGENTINA
AND
CMAT, UNIVERSIDAD DE LA REPÚBLICA
URUGUAY
E-MAIL: [fraiman@udesa.edu.ar](fraiman@udesa.edu.ar)